\documentclass[10pt]{amsart}
\usepackage[T1]{fontenc}
\usepackage{geometry}
\usepackage{amsmath}
\usepackage{amsfonts, amssymb, textcomp}
\usepackage[colorlinks=flase, linkcolor=red,urlcolor=green, citecolor=blue]{hyperref}
\usepackage{subeqnarray}

\usepackage{latexsym}
\usepackage{fancyhdr}
\usepackage{longtable}
\usepackage{amsmath, amssymb}
\usepackage{graphicx}

\numberwithin{equation}{section}

\newtheorem{lemma}[subsubsection]{Lemma}
\newtheorem{theorem}[subsubsection]{Theorem}
\newtheorem{proposition}[subsubsection]{Proposition}
\newtheorem{corollary}[subsubsection]{Corollary}
\newtheorem{definition}[subsubsection]{Definition}

\newtheorem{remark}[subsubsection]{Remark}

\newcommand{\RR}{\mathbb{R}}
\newcommand{\CC}{\mathbb{C}}
\newcommand{\NN}{\mathbb{N}}

\let\on=\operatorname

\newenvironment{demo}[1]{\par\smallskip\noindent{\bf #1.}}{\par\smallskip}

\title[How far is the Borel map from being surjective?]
{How far is the Borel map from being surjective in quasianalytic ultradifferentiable classes?}

\author[C.~Esser]{C\'eline Esser}

\author[G.~Schindl]{Gerhard Schindl}

\address{C.~Esser: Universit\'e de Li\`ege, D\'epartement de Math\'ematique, Quartier Polytech 1, All\'ee de la D\'ecouverte 12, B\^atiment B37, B-4000 Li\`ege, Belgique}
\email{celine.esser@uliege.be}

\address{G.~Schindl: Departamento de \'Algebra, An\'alisis Matem\'atico, Geometr{\'\i}a y Topolog{\'\i}a, Universidad de Valladolid, Facultad de Ciencias, Paseo de Bel\'en 7, 47011 Valladolid, Spain.}
\email{gerhard.schindl@univie.ac.at}

\begin{document}

\begin{abstract}
The Borel map $j^{\infty}$ takes germs at $0$ of smooth functions to the sequence of iterated partial derivatives at $0$. In the literature, it is well known that the restriction of $j^{\infty}$ to the germs of quasianalytic ultradifferentiable classes which are strictly containing the real analytic functions can never be onto the corresponding sequence space. In this paper, we are interested in studying how large the image of $j^{\infty}$ is and we investigate the size and the structure of this image by using different approaches (Baire residuality, prevalence and lineability). We give an answer to this question in the very general setting of quasianalytic ultradifferentiable classes defined by weight matrices, which contains as particular cases the classes defined by a single weight sequence or by a weight function.
\end{abstract}

\thanks{CE is supported by a F.R.S.-FNRS grant; GS is supported by FWF-Project J~3948-N35.}
\date{\today}

\maketitle

\section{Introduction}

In 1895, E. Borel proved that given any sequence $(a_n)_{n \in \NN}$ of complex numbers, there exists a infinitely differentiable function such that $f^{(n)}(0) = a_n$ for every $n \in \NN$ \cite{Borel}. This work has been investigated and extended ever since by many authors. In particular, the question has been handled in the context of so-called ultradifferentiable classes which are subclasses of smooth functions defined by imposing growth conditions on the derivatives of the functions using weight sequences $M$, functions $\omega$ or matrices $\mathcal{M}$, see~\cite{Carleman23, Carleman26, petzsche,thilliez,MR1150740,MR1021222,MR3073460,borelmappingquasianalytic}.\vspace{6pt}

Historically, those classes have been first introduced by using weight sequences, motivated among others by the characterization of the regularity of solutions of the heat equation or of other partial differential equations, see e.g. \cite{Rodino}. In order to measure the decay of the Fourier transform of smooth functions with compact support, classes of ultradifferentiable functions have then been defined using weight functions, e.g. see \cite{Beurling61} and \cite{PetzscheVogt}. In \cite{BraunMeiseTaylor90}, it turned out that such a behavior can also equivalently be expressed by having control on the growth of all the derivatives of the function itself in terms of this weight function and in \cite{BonetMeiseMelikhov07} it has been shown that classes defined in terms of weight sequences and weight functions are in general mutually distinct. Finally, in \cite{compositionpaper} and \cite{dissertation}, classes defined by weight matrices have been considered. It turned out that the weight sequence and weight function frameworks are particular cases of this setting, and this general method allows to treat both classical approaches jointly but also leads to more general classes.\vspace{6pt}

We say that an ultradifferentiable class is {\itshape quasianalytic} if the restriction of the Borel map $f \mapsto (\partial^\alpha f (0))_{\alpha \in \NN^r}$  to this class is injective; this notion plays an important role in many different contexts and applications (e.g. such classes do not contain partitions of unity). It came out of many studies that the restriction of the Borel map to the germs of quasianalytic ultradifferentiable classes which are strictly containing the real analytic functions can never be onto the corresponding sequence space. However, an interesting remaining question is ``how far away the Borel map is from being surjective?'' This is the question we tackle in this paper: We show that the image of the Borel map is ``small'' in the corresponding sequence space, using different approaches (as done e.g. in \cite{Esser14}). Let us present these different notions here.

\medskip

First, let us recall the following classical definition which gives a notion of residuality from a topological point of view.

\begin{definition}
If $X$ is a Baire space, then a subset $L \subset X$ is called {\em comeager} (or {\em residual}) if $L$ contains a countable intersection of dense open sets of $X$. The complement of a residual set is a {\em meager} (or {\em first category}) set in $X$.
\end{definition}

In order to get result about the ``size'' of sets from a measure-theoretical point of view, the notion of prevalence can be used. It has been introduced in \cite{Christensen74,HSY} to give an extension of the concept of ``almost everywhere'' (for the  Lebesgue measure) to metric infinite dimensional spaces (in these spaces,  no measure  is both  $\sigma$-finite  and translation invariant).

\begin{definition}\label{def:prev}
Let $X$ denote a complete metric vector space. A Borel subset $B\subset X$ is called \emph{Haar-null} if there exists a compactly supported probability measure $\mu$ such that
\begin{equation}\label{trans}
\forall x\in X,\quad \mu(x+B)=0.
\end{equation}
A subset $S$ of $X$ is called \emph{Haar-null} if it is contained in a Haar-null Borel set. A \emph{prevalent} set is the complement of a Haar-null set.
\end{definition}

The following results of \cite{Christensen74} and  \cite{HSY}  enumerate  important basic properties of prevalent sets:
\begin{itemize}
\item If $S$ is Haar-null, then $x+S$ is Haar-null for any $x \in X$.
\item If the dimension of $X$ is finite, $S$ is Haar-null if and only if $S$ has Lebesgue measure $0$.
\item Prevalent sets are dense.
\item Any countable intersection of prevalent sets is prevalent.
\end{itemize}

\begin{remark}
A useful way to get that a Borel set is Haar-null is to try the Lebesgue measure on the unit ball of a finite dimensional subspace $V$. In this context, condition (\ref{trans}) is equivalent to
\[ \forall x\in X, \;\;\;\; (x+B)\cap V \; \mbox{is of Lebesgue measure
zero}.\]
In this case, we say that $V$ is  a {\em probe} for the complement of $B$.
\end{remark}

Finally, we will also consider the notion of {\itshape lineability}, introduced in \cite{AGS}. This notion was motivated by the increasing interest toward the search for large algebraic structures of special objects (see \cite{BPS} for a review).

\begin{definition}
A set $L$ in a vector space $X$ is said to be \emph{lineable} in $X$ if $L \cup  \{0\}$ contains an infinite dimensional vector space. 
\end{definition}
Note that in the above definition, the considered vector space is generated by the \emph{finite} linear combinations of the elements of an infinite basis.

\medskip

The present paper is organized as follows: In Section \ref{germsultradifferentiableweightsequence}, we recall basic definitions and results concerning weight sequences $M$, the $M$-ultradifferentiable classes (of germs at $0$) and the associated sequence spaces which will be needed. In Section \ref{sequencecase}, we recall some important elements of the proof of \cite[Theorem 3]{thilliez} which gives the non-surjectivity of the Borel map in the quasianalytic setting, and we explain how to extend it in different directions. This Theorem allows us to obtain that the image of the Borel map is small (i.e. meager, Haar-null) in the Beurling case and its complement is lineable in both Roumieu and Beurling cases for quasianalytic classes strictly containing the real analytic functions. Let us mention that the above results are obtained in a more general context, since we actually prove that the image of the Borel map defined on {\itshape any} quasianalytic germ class associated with a sequence $M$ is small in any weighted sequence space associated with another quasianalytic weight sequence $N$ (assuming some mild standard assumptions on $M$ and $N$).\vspace{6pt}

Using the results and techniques developed in the single weight sequence case, we study in Section~\ref{matrixcase} the weight matrix case as well. Here the weight matrix is a one-parameter family of sequences $\mathcal{M}=\{M^{(\lambda)}: \lambda \in \RR_{>0}\}$, $M^{(\lambda)}$ again having some mild standard assumptions. All the results from the previous section are transferred to this more general setting with the same generality.\vspace{6pt}

Finally, in Section \ref{weighfunctioncase}, we treat the weight function case $\omega$. This is done by using the weight matrix setting and the fact that with each weight function $\omega$ (satisfying some mild standard growth conditions), we can associate a weight matrix $\Omega=\{W^{(\lambda)}: \lambda \in \RR_{>0}\}$ such that the ultradifferentiable classes defined by $\omega$ and $\Omega$ coincide (as locally convex vector spaces) and similarly for the corresponding sequence spaces, see \cite[Sect. 6, Sect. 7]{dissertation}, \cite[Sect. 5]{compositionpaper} and \cite{borelmappingquasianalytic}. Thus the results in this section can be seen as immediate Corollaries of the previous Section \ref{matrixcase}.\vspace{6pt}

Note that the presentation of this work and the standard assumptions on the weight structures are similar to the ones considered in \cite{borelmappingquasianalytic}. Moreover, throughout this paper, we write $\NN=\{0,1,\dots\}$, $\mathcal{E}(U)$ and $\mathcal{C}^{\omega}(U)$ shall denote respectively the class of all $\CC$-valued smooth functions and the class of all real analytic functions defined on non-empty open $U\subseteq\RR^r$.\vspace{6pt}

\section{Weight sequences and germs of ultradifferentiable functions}\label{germsultradifferentiableweightsequence}

\subsection{Denjoy-Carleman  ultradifferentiable classes and their germs}

\begin{definition}
Let $M=(M_p)_{p \in \NN}\in\RR_{>0}^{\NN}$ be an arbitrary sequence of positive real numbers. Let $r\in\NN_{>0}$ and $U\subseteq\RR^r$ be non-empty and open. The $M$-\emph{ultradifferentiable Roumieu type class} is defined by
\begin{equation*}\label{roumieu}
\mathcal{E}_{\{M\}}(U):=\{f\in\mathcal{E}(U):\;\forall\;K\subseteq U\;\text{compact}\;\; \exists\;h>0,\;\|f\|^M_{K,h}<+\infty\},
\end{equation*}
and the $M$-\emph{ultradifferentiable Beurling type class} by
\begin{equation*}\label{beurling}
\mathcal{E}_{(M)}(U):=\{f\in\mathcal{E}(U):\;\forall\;K\subseteq U\;\text{compact}\;\; \forall\;h>0,\;\|f\|^M_{K,h}<+\infty\},
\end{equation*}
where (using the standard multi-index notation for the partial derivatives)
\begin{equation*}\label{semi-norm-1}
\|f\|^M_{K,h}:=\sup_{\alpha\in\NN^r,x\in K}\frac{|\partial^{\alpha}f(x)|}{h^{|\alpha|} M_{|\alpha|}}.
\end{equation*}
\end{definition}

As usual, we  will write $m=(m_p)_{p \in \NN}$ for $m_p:=\frac{M_p}{p!}$.

\begin{remark}
At this point, we want to make the reader aware that the sequence $M$ considered in  \cite{thilliez} is precisely the sequence $m=(m_p)_{p\in \NN}$ in the notation of this work.
\end{remark}

For any compact set $K$ with smooth boundary $\mathcal{E}_{M,h}(K):=\{f\in\mathcal{E}(K): \|f\|^M_{K,h}<+\infty\}$ is a Banach space.
The Roumieu type class is endowed with the projective topology w.r.t. all $K\subseteq U$ compact and the inductive topology w.r.t. $h\in\NN_{>0}$, whereas the Beurling type class is endowed with the projective topology w.r.t. $K\subseteq U$ compact and w.r.t. $1/h$, $h\in\NN_{>0}$. Hence $\mathcal{E}_{(M)}(U)$ is a {\itshape Fr\'echet space} and $\underset{h>0}{\varinjlim}\;\mathcal{E}_{M,h}(K)=\underset{n\in\NN_{>0}}{\varinjlim}\;\;\mathcal{E}_{M,n}(K)$ is a {\itshape Silva space}, i.e. a countable inductive limit of Banach spaces with compact connecting mappings, see \cite[Proposition 2.2]{Komatsu73}.\vspace{6pt}

Note that the special case $M_p=p!$ yields $\mathcal{E}_{\{M\}}(U)=\mathcal{C}^{\omega}(U)$, whereas $\mathcal{E}_{(M)}(U)$ consists of the restrictions of all entire functions provided that $U$ is connected.

\begin{definition}
The spaces of \emph{germs at $0\in\RR^r$} of the $M$-ultradifferentiable functions of Roumieu and Beurling types are defined respectively by
\begin{equation*}\label{roumieugerm}
\mathcal{E}_{\{M\}}^{0,r}:=\underset{k\in\NN_{>0}}{\varinjlim}\mathcal{E}_{\{M\}}\left(\left(-\frac{1}{k},\frac{1}{k}\right)^r\right),
\end{equation*}
and
\begin{equation*}\label{beurlinggerm}
\mathcal{E}_{(M)}^{0,r}:=\underset{k\in\NN_{>0}}{\varinjlim}\mathcal{E}_{(M)}\left(\left(-\frac{1}{k},\frac{1}{k}\right)^r\right).
\end{equation*}
\end{definition}

Again, if one considers the sequence $M_p=p!$ in the Roumieu case, we obtain the space of {\itshape germs of real analytic functions at} $0\in\RR^r$; it is denoted by $\mathcal{O}^{0,r}$.

\medskip
Let us now introduce the corresponding spaces of  complex sequences.

\begin{definition}
We define the sequence spaces $\Lambda^r_{\{ M\} }$ and $\Lambda^r_{(M) }$ by setting
$$
\Lambda^r_{\{ M\} }=\left\{\mathbf{b}\in\CC^{\NN^r} : \exists h > 0 , |\mathbf{b}|^M_{h}< + \infty \right\}
$$
and
$$
\Lambda^r_{( M) }=\left\{ \mathbf{b}\in\CC^{\NN^r} : \forall h > 0 , |\mathbf{b}|^M_{h}< + \infty \right\},
$$
where for any $h >0$,
$$
|\mathbf{b}|^M_{h} := \sup_{\alpha\in \NN^r}\frac{|b_{\alpha}|}{h^{|\alpha|} m_{|\alpha|}}=\sup_{\alpha\in \NN^r}\frac{|b_{\alpha}||\alpha|!}{h^{|\alpha|} M_{|\alpha|}}.
$$
\end{definition}

We endow these spaces also with their natural topology: $\Lambda_{\{M\} }$ is an (LB)-space and $\Lambda_{( M) }$ a Fr\'echet space.

\begin{remark}
Note that there exists a one-to-one correspondence between $\Lambda^r_{[ M] }$ and the (ring) of weighted formal power series, whose elements $F=\sum_{\alpha\in\NN^r}F_{\alpha}x^{\alpha}$ satisfy $|F_{\alpha}|\le Ch^{|\alpha|}m_{|\alpha|}$ for some $C,h>0$ and all $\alpha\in\NN^r$ in the Roumieu case, resp. for all $h>0$ small, some $C=C_h$ large and all $\alpha\in\NN^r$ in the Beurling case (as it as been considered in \cite[Section 1.2]{thilliez} for the Roumieu case).
\end{remark}

Finally, let us define the Borel map. For reasons of convenience, the following {\itshape convention} will also be used: we write $\mathcal{E}_{[M]}$ if either $\mathcal{E}_{\{M\}}$ or $\mathcal{E}_{(M)}$ is considered, but not mixing the cases if statements involve more than one $\mathcal{E}_{[M]}$ symbol. We use similar notations for the sequence classes $\Lambda^r_{[M]}$. In both cases, the \emph{Borel map} $j^\infty$ is defined by
\begin{equation}\label{Borelmapsequence}
j^{\infty}:\mathcal{E}_{[M]}^{0,r}\longrightarrow\Lambda^r_{[M]},\hspace{20pt}j^\infty(f) = \left(\frac{\partial^{\alpha}f(0)}{|\alpha|!}\right)_{\alpha\in \NN^r}.
\end{equation}

\begin{remark} At this step, let us point out that our definitions of sequence spaces differ from those in \cite[Section 2.1]{borelmappingquasianalytic}; this is due to technical reasons since we will have to work mainly with the sequence $m=(\frac{M_p}{p!})_{p \in \NN}$ instead of $M=(M_p)_{p\in\NN}$.
\end{remark}

\begin{remark}
Moreover, we outline that we could replace in the definition of the germ spaces the point $0\in\RR^r$ by any other point $a\in\RR^r$ and define the corresponding Borel map $j^{\infty}_a$ similarly as in \eqref{Borelmapsequence}.
\end{remark}

\subsection{Weight sequences}\label{weightsequences}

We consider the following definition, according to \cite[Section 2.2]{borelmappingquasianalytic}.

\begin{definition}\label{def_weightseq}
A sequence of positive real numbers $M=(M_p)_{p \in \NN} \in\RR_{>0}^{\NN}$ is called a \emph{weight sequence} if
\begin{itemize}
\item[$(I)$] $1=M_0\le M_1$ (normalization),
\item[$(II)$] $p\mapsto M_p$ is log-convex, or equivalently $M_p^2\le M_{p-1}M_{p+1}$ for all $p\in \NN_{>0}$,
\item[$(III)$] $\liminf_{p\rightarrow\infty}(m_p)^{1/p}>0$.
\end{itemize}
Recall that $m_p:=\frac{M_p}{p!}$ for every $p \in \NN$.
\end{definition}

Let us recall that if $M$ is log-convex and normalized, then $M$ and $k\mapsto(M_k)^{1/k}$ are both increasing and $M_jM_k\le M_{j+k}$ holds for all $j,k\in\NN$, e.g. see \cite[Lemmata 2.0.4, 2.0.6]{diploma}.

\medskip

Given two (weight) sequences, we write $M\le N$ if and only if $M_p\le N_p$ holds for all $p\in\NN$ and define the relations
$$M\hypertarget{mpreceq}{\preceq}N:\Leftrightarrow\;\exists\;C\ge 1\;\forall \, p\in\NN:\; M_p\le C^pN_p\Longleftrightarrow\sup_{p\in\NN_{>0}}\left(\frac{M_p}{N_p}\right)^{1/p}<+\infty,$$
$$M\hypertarget{mtriangle}{\vartriangleleft}N:\Leftrightarrow\;\forall\;h>0\;\exists\;C_h\ge 1\;\forall \, p\in\NN:\; M_p\le C_h h^p N_p\Longleftrightarrow\lim_{p\rightarrow\infty}\left(\frac{M_p}{N_p}\right)^{1/p}=0.$$
It is straightforward to see that in the above relations we can replace the sequences $M$ and $N$ simultaneously by the sequences $m$ and $n$.

\medskip

Those relations between weight sequences imply inclusions between ultradifferentiable classes, e.g. see \cite[Section 2.2]{borelmappingquasianalytic} and the references therein.

More precisely, let $M$ be a weight sequence and $N$ arbitrary, then $M\hyperlink{mpreceq}{\preceq}N$ if and only if $\mathcal{E}_{[M]}\subseteq\mathcal{E}_{[N]}$, which is equivalent to $\Lambda^r_{[M]}\subseteq\Lambda^r_{[N]}$. In particular, choosing $M=(p!)_{p\in\NN}$, we get $\mathcal{C}^{\omega}\subseteq\mathcal{E}_{\{N\}}$ if and only if $\liminf_{p\rightarrow + \infty}(n_p)^{1/p}>0$. Moreover, if $N$ is a weight sequence, then $\mathcal{E}_{\{N\}}\subseteq\mathcal{C}^{\omega}$ if and only if $\sup_{p\in\NN_{>0}}(n_p)^{1/p}<+\infty$. Hence $\mathcal{C}^{\omega}\subsetneq\mathcal{E}_{\{N\}}$ if and only if $\sup_{p\in\NN_{>0}}(n_p)^{1/p}=+\infty$.

Similarly $M\hyperlink{mtriangle}{\vartriangleleft}N$ if and only if $\mathcal{E}_{\{M\}}\subsetneq\mathcal{E}_{(N)}$, which is equivalent to $\Lambda^r_{\{M\}}\subsetneq\Lambda^r_{(N)}$. In particular, $\mathcal{C}^{\omega}\subsetneq\mathcal{E}_{(N)}$ if and only if $\lim_{p\rightarrow + \infty}(n_p)^{1/p}=+\infty$.

\begin{definition}
A weight sequence $M$ is called {\itshape quasianalytic} if
\begin{equation}\label{eq_Q}
\sum_{p=1}^{+ \infty}\frac{M_{p-1}}{M_p}=+\infty . \tag{Q}
\end{equation}
\end{definition}

By using {\itshape Carleman's inequality} (a proof is presented in \cite[Proposition 4.1.7]{diploma}), one can show that
$$\sum_{p=1}^{+ \infty}\frac{M_{p-1}}{M_p}=+\infty\Longleftrightarrow\sum_{p=1}^{+ \infty}\frac{1}{(M_p)^{1/p}}=+\infty.$$

\begin{definition}
A subclass $\mathcal{Q}\subseteq\mathcal{E}$ is called {\itshape quasianalytic} if for any open connected set $U\subseteq\RR^r$ and each point $a\in U$, the Borel map $j^{\infty}_a$ is {\itshape injective} on $\mathcal{Q}(U)$.
\end{definition}

In the case $\mathcal{Q}\equiv\mathcal{E}_{[M]}$ the {\itshape Denjoy-Carleman theorem} characterizes this behavior in terms of the defining weight sequence $M$. More precisely, it states that $\mathcal{E}_{[M]}$ is quasianalytic if and only if $M$ satisfies (\ref{eq_Q}). Let us moreover mention  that $\mathcal{E}_{[M]}$ is quasianalytic if and only if there do not exist non-trivial functions in $\mathcal{E}_{[M]}$ with compact support, e.g. see \cite[Thm. 19.10]{rudin}.

A basic assumption in the proof of the Denjoy-Carleman theorem is condition $(II)$ for $M$ (i.e. $M$ is log-convex). But sometimes it might be convenient to skip this condition and to work in a more general stetting: In this case, one considers admissible regularizations of $M$ (with $M_0=1$), see \cite[Section 4.3]{testfunctioncharacterization} and the references therein. For this reason, we shall denote by
\begin{itemize}
\item $M^{\on{lc}}$ the log-convex minorant of $M$, i.e. the largest sequence $N$ such that $N$ is log-convex and $N\le M$,
\item $M^I$ the sequence for which $((M^I_p)^{1/p})_{p \in \NN_{>0}}$ is the increasing minorant of $((M_p)^{1/p})_{p \in \NN_{>0}}$ (and put $M^I_0:=1$).
\end{itemize}
Consequently, $M=M^{\on{lc}}$ if and only if $(II)$, and $M^I=M$ if and only if $k\mapsto(M_k)^{1/k}$ is increasing. One has $M^{\on{lc}}\le M^I\le M$ since a log-convex weight sequence is increasing. We can now recall the following result, see \cite[Proposition 4.4]{testfunctioncharacterization} which is based on \cite[Theorem 1.3.8]{hoermander} for the Roumieu case, and see \cite[Theorem 4.2]{Komatsu73} for the Beurling case.

\begin{proposition}\label{nonquasiremarks}
Let $M\in\RR_{>0}^{\NN}$ with $M_0=1$. The following assertions are equivalent:
\begin{itemize}
\item $\mathcal{E}_{[M]}$ is quasianalytic,
\item $M^{\on{lc}}$ satisfies (\ref{eq_Q}),
\item $\sum_{p= 1}^{+ \infty}\frac{1}{(M^I_p)^{1/p}}=+\infty$.
\end{itemize}
\end{proposition}

\begin{remark}\label{logminorantremark}
We mention that in the following sections, we will study the Borel map $j^{\infty}$ defined in quasianalytic ultradifferentiable classes such that $\mathcal{C}^{\omega}\subsetneq\mathcal{E}_{[M]}$ holds true. As already pointed out in \cite[Remark 1]{borelmappingquasianalytic}, the general assumptions $(I)-(III)$ on $M$ are not restricting the generality of our considerations. More precisely, for any $M\in\RR_{>0}^{\NN}$ with $\mathcal{C}^{\omega}\subseteq\mathcal{E}_{[M]}$ we have $\liminf_{p\rightarrow + \infty}(m_p)^{1/p}>0$ in the Roumieu and $\lim_{p\rightarrow+ \infty}(m_p)^{1/p}=+\infty$ in the Beurling case (see also \cite[Prop. 2.12 $(4),(5)$]{compositionpaper}). Then, by \cite[Theorem 2.15]{compositionpaper}, we can replace $M$ by $M^{\on{lc}}$ without changing the associated ultradifferentiable class whereas only $\Lambda^r_{[M^{\on{lc}}]}\subseteq\Lambda^r_{[M]}$ follows (and the weight matrix/function setting is reduced to the sequence case situation as will be seen in the next sections).
\end{remark}

\begin{remark}
Let us point out that all results below also hold true if $0\in\RR^r$ is replaced by any other point $a\in\RR^r$ (translation).
\end{remark}

\section{The weight sequence case $M$}\label{sequencecase}

\subsection{Thilliez's proof for non-surjectivity}\label{generalizationThilliezsequence}
Let $M$ be a weight sequence, i.e. satisfying our standard assumptions $(I)-(III)$. To ensure that the real analytic functions/germs are strictly contained in the considered class $\mathcal{E}_{[M]}$, we have to assume
$$\sup_{k\in\NN_{>0}}(m_k)^{1/k}=+\infty\hspace{30pt}\text{or}\hspace{30pt}\lim_{k\rightarrow + \infty}(m_k)^{1/k}=+\infty$$
in the Roumieu or in the Beurling case respectively.

\medskip

The aim of this subsection is to recall the main elements of the proof of \cite[Theorem 3]{thilliez}, which is based on the original ideas of Carleman \cite{Carleman26}, to be applicable in our present context. We will also explain how it can be extended: Indeed, in \cite{thilliez} only the Roumieu case has been treated and it has been assumed there that $m=(m_k)_{k \in \NN}$ is log-convex (in that case we say that $M$ is strongly log-convex), which implies that $k\mapsto(m_k)^{1/k}$ is increasing. Consequently, also in the Roumieu case, the assumption for the strict inclusion turns into $\lim_{k\rightarrow + \infty}(m_k)^{1/k}=+\infty$.

\medskip

In our approach we do not want to assume strongly log-convexity on $M$, or more generally on some or all $M^{(\lambda)}\in\mathcal{M}$ in the weight matrix case considered in Section \ref{matrixcase} below: This is due to the fact that on the one hand, in general we do not know whether some or all of the sequences $W^{(\lambda)}$ of the matrix $\Omega$ associated with a weight function $\omega$ will satisfy this requirement, see Section \ref{weighfunctioncase} below for further explanations. On the other hand, in any cases, strongly log-convexity seems to be too strong and superfluous in studying the questions under consideration in this paper (but not for some questions studied in \cite{thilliez}).

A second generalization is that we also consider a kind of mixed setting of two (in general different) weight sequences $M$ and $N$.

\vspace{6pt}

Let us start by recalling the following representation formula, obtained within the first part of the proof of \cite[Theorem 3]{thilliez}. As mentioned before, this result has been obtained by assuming the strongly log-convexity on $M$. However, by following directly the lines of this proof, the result still holds with the weaker (basic) assumptions on $M$. Hence, we have the following Theorem.

\begin{theorem}[Representation formula, \cite{thilliez}]\label{prop_representationformula}
Let $M$ be a quasianalytic weight sequence. There exist numbers $(\omega^M_{j,k})_{j,k \in \NN}$ such that
\begin{equation}\label{eq:lim_omega}
\lim_{k \to + \infty} \omega^M_{j,k} =1,  \quad \forall j \in \NN,
\end{equation}
and such that, given any function $f \in \mathcal{E}^{0,1}_{\{M\}}$, one has
\begin{equation}\label{representation formula}
f(x) = \lim_{k \to + \infty} \sum_{j=0}^{k-1} \omega^M_{j,k} \frac{f^{(j)} (0)}{j !} x^j
\end{equation}
for every $x>0$ small enough.
\end{theorem}

Keeping the notations of this Theorem, we directly get the following important result.

\begin{corollary}\label{cor_Thilliez}
Let $M$ be a quasianalytic weight sequence. If $\mathbf{b}=(b_j)_{j \in \NN} \in \CC^{\NN}$ is a sequence for which there exists a sequence of positive real numbers $(a_n)_{n \in \NN}$ decreasing to $0$ such that
\begin{equation}\label{eq_cor_Thilliez}
\limsup_{k \to + \infty} \left| \sum_{j=0}^{k-1} \omega^{M}_{j,k} b_j a_n^j \right|= + \infty
\end{equation}
for all $n \in \NN$, then $\mathbf{b} \notin j^\infty (\mathcal{E}^{0,1}_{\{M\}})$.
\end{corollary}

\demo{Proof}
Assume by contradiction that we can find $f\in\mathcal{E}_{\{M\}}^{0,1}$  such that $j^{\infty}(f)=\mathbf{b}$. Using the representation formula \eqref{representation formula} of Theorem \ref{prop_representationformula} together with the definition of the Borel map, we get
\begin{equation*}\label{omega2}
f(x)=\lim_{k\rightarrow + \infty}\sum_{j=0}^{k-1}\omega^M_{j,k}b_jx^j
\end{equation*}
for every $x>0$ small enough, hence a contradiction.
\qed\enddemo

\begin{remark}\label{rem:dim_r}
Let $\mathbf{b}=(b_j)_{j \in \NN}$ be a sequence which does not belong to $j^\infty (\mathcal{E}^{0,1}_{\{M\}})$. Then, for $r>1$, any sequence $\mathbf{\widetilde{b}} \in \CC^{\NN^r}$ satisfying $\widetilde{b}_{(j,0,\dots,0)} = b_j$ for every $j \in \NN$  is not in the image $j^\infty (\mathcal{E}^{0,r}_{\{M\}})$: Indeed, if one assumes now that there is $f \in \mathcal{E}^{0,r}_{\{M\}}$ such that $\mathbf{\widetilde{b}} = j^{\infty}(f)$, then by considering the restriction mapping $R:\mathcal{E}_{\{M\}}^{0,r}\twoheadrightarrow\mathcal{E}_{\{M\}}^{0,1}$, the restriction $R(f)$ of $f$  would belong to $\mathcal{E}_{\{M\}}^{0,1}$ and one would obtain that $j^{\infty}(R(f)) = \mathbf{b}$.

\end{remark}

The following Theorem is a direct generalization  of the second part of the proof of \cite[Theorem 3]{thilliez}: We consider two weight sequences (different or not) and we treat the Beurling case as well.

\begin{theorem}\label{Thillieztheorem}
Let $M$ and $N$ be two quasianalytic weight sequences such that $\sup_{k\in\NN_{>0}}(n_k)^{1/k}=+\infty$ resp. $\lim_{k\rightarrow + \infty}(n_k)^{1/k}=+\infty$, i.e. $\mathcal{O}^{0,r}\subsetneq\mathcal{E}_{\{N\}}^{0,r}$ resp. $\mathcal{O}^{0,r}\subsetneq\mathcal{E}_{(N)}^{0,r}$. Then, one has $$j^{\infty}(\mathcal{E}_{\{M\}}^{0,r}) \cap \Lambda^r_{[N]} \subsetneq \Lambda^r_{[N]}$$ (and hence $j^{\infty}(\mathcal{E}_{(M)}^{0,r}) \cap \Lambda^r_{[N]} \subsetneq \Lambda^r_{[N]}$ also). In particular, the Borel map $j^{\infty}:\mathcal{E}_{[N]}^{0,r}\longrightarrow\Lambda^r_{[N]}$ is not surjective.
\end{theorem}

\begin{remark} Theorem \ref{Thillieztheorem} is stronger than only having non-surjectivity of $j^{\infty}:\mathcal{E}_{[N]}^{0,r}\longrightarrow\Lambda^r_{[N]}$ since $M$ can be any other quasianalytic weight sequence satisfying $N\hyperlink{mpreceq}{\preceq}M$ (i.e. much larger than $N$).
\end{remark}

This theorem will follow directly from the Corollary \ref{cor_Thilliez} and the next lemma which gives the existence of sequences satisfying (\ref{eq_cor_Thilliez}) in any class $\Lambda^r_{[N]}$. The proof of this lemma reduces to the argument given in the proof of \cite[Theorem 3]{thilliez} with the only difference that the sequence $M$ is replaced by the square root of the sequence $(n_j)_{j \in \NN}$ to treat both the Beurling case and the mixed setting.

\begin{lemma}\label{lem_Thilliez}
Let $M$ and $N$ be two quasianalytic weight sequences such that $\sup_{k\in\NN_{>0}}(n_k)^{1/k}=+\infty$ resp. $\lim_{k\rightarrow + \infty}(n_k)^{1/k}=+\infty$, i.e. $\mathcal{O}^{0,1}\subsetneq\mathcal{E}_{\{N\}}^{0,1}$ resp. $\mathcal{O}^{0,1}\subsetneq\mathcal{E}_{(N)}^{0,1}$. There exists $\mathbf{F} \in \Lambda^1_{[N]}$ and $a_0  \in (0,1]$  such that
\begin{equation*}\label{eq_cor_Thilliez1}
\limsup_{k \to + \infty} \left| \sum_{j=0}^{k-1} \omega^{M}_{j,k} F_j a^j \right|= + \infty,
\end{equation*}
for all $0<a \leq a_0$.
\end{lemma}

Since it will be useful in the next section, let us recall that such a sequence $\mathbf{F}$ can be obtained by setting
\begin{equation}\label{squenceF}
F_{k_p}:= \sqrt{n_{k_p}} \quad \text{ and } \quad F_j := 0 \, \text{ otherwise,}
\end{equation}
where the increasing sequence $(k_p)_{p \in \NN}$ of natural numbers is chosen such that
\begin{equation}\label{powerseriesrad}
\sup_{p\in\NN}\left(n_{k_p}\right)^\frac{1}{k_p}=+\infty\hspace{10pt}\text{resp.}\hspace{10pt}\lim_{p\rightarrow + \infty}\left(n_{k_p}\right)^\frac{1}{k_p}=+\infty
\end{equation}
and
\begin{equation}\label{sequencelequ2}
\sum_{j=0}^{k_{p-1}}\left|\omega^M_{j,k_p}-1\right||F_j| \le 1, \quad \forall p \in \NN_{>0}.
\end{equation}

\begin{remark}
As defined in \eqref{squenceF}, the sequence $\mathbf{F}=(F_j)_{j \in \NN}$ can never define a real analytic germ, otherwise $F_{k_p}=(n_{k_p})^{1/2}\le C h^{k_p}$ should be satisfied for some $C,h>0$ and all $p\in\NN$, a contradiction to \eqref{powerseriesrad} in the Roumieu and the Beurling case.
\end{remark}

\begin{remark}\label{remarklpj}
Let us note that it is possible to define a sequence $\mathbf{F}$ satisfying the assumption of Lemma \ref{lem_Thilliez} with only non-zero elements by setting  $F_j := l^p_j$ for $k_{p-1}<j<k_p$, where the values $l^p_j$ are subjected to some precise growth control. This provides some additional information on sequences not contained in the image of the Borel map.

In particular, negative values $l^p_j$ or mixed signs are allowed : we meet here a situation not treated in \cite{borelmappingquasianalytic}.

However, the construction of $\mathbf{F}$ depends on the given $M$. Even if $l^p_j\equiv 0$ for all $p,j$, as in the original proof of \cite[Theorem 3]{thilliez} (or as in Lemma \ref{lem_Thilliez} above), then \eqref{sequencelequ2} and hence the choice of $(k_p)_{p \in \NN}$ is still depending on $M$. It is not clear how to get rid of this problem in general. 
Consequently it seems not possible to prove by using this technique the existence of $\mathbf{F}=(F_j)_{j \in \NN}\in \Lambda^1_{[N]}$ such that $\mathbf{F}$ does not belong to $j^\infty (\mathcal{E}_{\{M\}}^{0,1})$ for {\itshape any} quasianalytic weight sequence $M$ as it has been done in \cite[Thm. 2, Thm. 3]{borelmappingquasianalytic} by using completely different methods (and sequences $\mathbf{F}=(F_j)_{j \in \NN}$ such that $F_j>0$ for all $j\in\NN$ and not defining a real analytic germ).

Since this more general definition of $\mathbf{F}$ would neither change nor simplify the proofs of the main results below and would unnecessarily complicate the notation we will work with lacunary sequences as stated above, i.e. $l^p_j\equiv 0$ for all $p,j$.
\end{remark}

\subsection{Generic size of the image of the Borel map}\label{Borelmapfaraway}

Let $M$ and $N$ be two (in general different) quasianalytic weight sequences. The aim of this section is to study the size of $j^\infty (\mathcal{E}^{0,r}_{[M]})$ in $\Lambda^r_{[N]}$ using the different notions of genericity presented in the introduction. The results will be obtained by applying Corollary \ref{cor_Thilliez}, Remark \ref{rem:dim_r} and Lemma \ref{lem_Thilliez}.\vspace{6pt}

First, let us concentrate on the Beurling case. Indeed, we intend to study the size of the image from the point of view of {\itshape Baire genericity} (resp. {\itshape prevalence}), for which the underlying space needs to be a Baire space (resp. a complete metrizable space). In the next Theorem, we prove that $\Lambda^r_{(N)} \setminus j^\infty (\mathcal{E}^{0,r}_{\{M\}})$ (and hence also $\Lambda^r_{(N)} \setminus j^\infty (\mathcal{E}^{0,r}_{(M)})$) is a ``big'' set in  $\Lambda^r_{(N)}$, hence the image of the Borel mapping defined on {\itshape any} quasianalytic class of $M$-ultradifferentiable germs is ``small'' in the space $\Lambda^r_{(N)}$.

\begin{theorem}\label{thm:Baire}
Let $M$ and $N$ be two quasianalytic weight sequences.

Let us assume that $\lim_{k\rightarrow + \infty}(n_k)^{1/k}=+\infty$, i.e. $\mathcal{O}^{0,r}\subsetneq\mathcal{E}^{0,r}_{(N)}$. Then, the set $j^{\infty}(\mathcal{E}_{\{M\}}^{0,r}) \cap \Lambda^r_{(N)} $ is meager in $\Lambda^r_{(N)}$, i.e. $\Lambda^r_{(N)} \setminus j^\infty (\mathcal{E}^{0,r}_{\{M\}})$ contains a countable intersection of dense open sets.
\end{theorem}

Note that as a particular case, the choice $M\equiv N$ yields that $\Lambda^r_{(M)} \setminus j^\infty (\mathcal{E}^{0,r}_{\{M\}})$ contains a countable intersection of dense open sets.

\demo{Proof}
By Lemma \ref{lem_Thilliez}, we can consider $\mathbf{F} \in \Lambda^r_{(N)}$ and $a_0 \in (0,1]$ such that
\begin{equation}\label{thm:Baire:eq30}
\limsup_{k \to + \infty} \left| \sum_{j=0}^{k-1} \omega^{M}_{j,k} F_j a^j \right| = + \infty
\end{equation}
for any $0<a \leq a_0$, where $F_j = F_{(j,0,\dots,0)}$ for any $j \in \NN$. Let us fix a sequence $ (a_p)_{p \in \NN}$ in $(0,a_0]$ decreasing to $0$. If the set $\mathcal{G}$ is defined by
\begin{equation*}\label{eq:Ga2}
\mathcal{G} =\bigcap_{p \in \NN} \left\{\mathbf{b} \in \Lambda^r_{(N)} : \limsup_{k \to + \infty} \left| \sum_{j=0}^{k-1} \omega^{M}_{j,k} b_j a_p^j \right|= + \infty  \right\} ,
\end{equation*}
we know from Corollary \ref{cor_Thilliez} and Remark \ref{rem:dim_r} that
\begin{equation}\label{eq:Ga1}
\mathcal{G} \subseteq \Lambda^r_{(N)} \setminus j^\infty (\mathcal{E}^{0,r}_{\{M\}}) \, .
\end{equation}

In order to get the result, it suffices then to prove that $\mathcal{G}$ can be written as a countable intersection of dense open sets of $\Lambda^r_{(N)}$. One has
\begin{equation*}
\mathcal{G} = \bigcap_{p \in \NN} \bigcap_{P \in \NN} \bigcap_{K \in \NN_{>0}} \bigcup_{k \ge K} G(p,P,k)
\end{equation*}
where
\begin{equation*}
G(p,P,k) =  \left\{\mathbf{b} \in \Lambda^r_{(N)} : \left| \sum_{j=0}^{k-1} \omega^{M}_{j,k} b_j a_p^j \right| > P  \right\}.
\end{equation*}
Let us fix $p \in \NN$, $P \in \NN$ and $k \in \NN_{>0}$ and let us show that $G(p,P,k)$ is open. Let us consider a sequence  $(\mathbf{b}^{(l)})_{l\in \NN}$ of elements of $\Lambda^r_{(N)}$ which does not belong to $G(p,P,k) $ and which converges to $\mathbf{b}$ in $\Lambda^r_{(N)}$, and let us show that
\begin{equation}\label{thm:Baire:eq1}
\left| \sum_{j=0}^{k-1} \omega^{M}_{j,k} b_j a_p^j \right| \leq P.
\end{equation}
Let $\delta>0$ be arbitrary and fix $\varepsilon >0$ such that
\begin{equation*}
\varepsilon < \frac{\delta}{\sum_{j=0}^{k-1} | \omega^{M}_{j,k}| n_j a_p^j  }.
\end{equation*}
By assumption, there is $L \in \NN$ such that $\big| \mathbf{b}^{(l)} - \mathbf{b} \big|^{N}_1 \leq \varepsilon$  for all $l \geq L$. Then, for all $l \geq L$, one has
\begin{align*}
 \left| \sum_{j=0}^{k-1} \omega^{M}_{j,k} b_j a_p^j \right|  & \leq  \left| \sum_{j=0}^{k-1} \omega^{M}_{j,k} (b_j^{(l)} -b_j) a_p^j \right|  +  \left| \sum_{j=0}^{k-1} \omega^{M}_{j,k} b^{(l)}_j a_p^j \right| \\
  & \leq   \sum_{j=0}^{k-1} | \omega^{M}_{j,k} | \big|b_j^{(l)} -b_j \big| a_p^j  + P  \\
   & \leq  \varepsilon \sum_{j=0}^{k-1} | \omega^{M}_{j,k} | n_j a_p^j  + P   \leq \delta + P ,
 \end{align*}
hence (\ref{thm:Baire:eq1}) since $\delta >0$ is arbitrary. It follows that for any $p \in \NN$, $P \in \NN$ and $K \in \NN_{>0}$
$$
\bigcup_{k \ge K} G(p,P,k)
$$
is open, and it remains to prove that it is dense in $\Lambda^r_{(N)}$. So, let us consider an arbitrary $\mathbf{b} \in \Lambda^r_{(N)}$ and let us fix $\varepsilon >0$. It follows from (\ref{thm:Baire:eq30}) that for all $K \in \NN$, there is $k \geq K$ such that
\begin{equation}\label{thm:Baire:eq2}
\left| \sum_{j=0}^{k-1} \omega^{M}_{j,k} F_j a_p^j \right| \geq \frac{P}{\varepsilon}.
\end{equation}
Then, either $\mathbf{b} + \varepsilon \mathbf{F}$ or $\mathbf{b} - \varepsilon \mathbf{F}$ belongs to $\bigcup_{k \ge K} G(p,P,k)$: Otherwise, one would have
$$
2 \varepsilon \left| \sum_{j=0}^{k-1} \omega^{M}_{j,k} F_j a_p^j \right|
\leq   \left| \sum_{j=0}^{k-1} \omega^{M}_{j,k} (b_j-  \varepsilon F_j) a_p^j \right| +  \left| \sum_{j=0}^{k-1} \omega^{M}_{j,k} ( b_j+ \varepsilon F_j  ) a_p^j \right|  \leq 2P
$$
which contradicts (\ref{thm:Baire:eq2}). Moreover, for any $h>0$,  one has
\begin{equation}\label{thm:Baire:eq3}
| \mathbf{b} - (\mathbf{b} \pm \varepsilon \mathbf{F})|^{N}_h = \varepsilon |\mathbf{F}|^{N}_h,
\end{equation}
and the density of $\bigcup_{k \ge K} G(p,P,k)$ in $\Lambda^r_{(N)}$ follows.
\qed\enddemo

\begin{remark}\label{rem:Baire}
In particular, from Baire's theorem, one gets that the set $\Lambda^r_{(N)} \setminus j^\infty (\mathcal{E}^{0,r}_{\{M\}})$ is dense in $\Lambda^r_{(N)}$.
\end{remark}

Let us now show that the previous result also holds in the context of prevalence.

\begin{theorem}\label{thm:prevalence}
Let $M$ and $N$ be two quasianalytic weight sequences.

Let us assume that $\lim_{k\rightarrow + \infty}(n_k)^{1/k}=+\infty$, i.e. $\mathcal{O}^{0,r}\subsetneq\mathcal{E}^{0,r}_{(N)}$. Then, the set $j^{\infty}(\mathcal{E}_{\{M\}}^{0,r}) \cap \Lambda^r_{(N)} $ is Haar-null in $\Lambda^r_{(N)}$.
\end{theorem}

\demo{Proof}
We use similar notations as in the proof of Theorem \ref{thm:Baire}. From (\ref{eq:Ga1}), it suffices to prove that $\mathcal{G}$ is prevalent in $\Lambda^r_{(N)}$. We already know that it is a Borel set, since it is a countable intersection of the open sets $ \bigcup_{k \ge K} G(p,P,k)$. Let us prove that each of these sets is prevalent, hence the result since a countable intersection of prevalent sets is prevalent. We use for a probe the space generated by $\mathbf{F}$. For any $\mathbf{b} \in \Lambda^r_{(N)}$, the line $L := \big\{\mathbf{b}+ \alpha \mathbf{F} : \alpha \in \RR \big\}$ contains at most one element in the set $\Lambda^r_{(N)} \setminus \bigcup_{k \ge K} G(p,P,k)$. Indeed, assume that there exist two different such sequences in $L$ associated with the numbers $\alpha,\beta \in \RR$, $\alpha\neq\beta$. Then, for all $k \geq K$, one would have
$$
\left| \sum_{j=0}^{k-1} \omega^{M}_{j,k} (b_j + \alpha F_j) a_p^j \right|  \leq P \quad \text{and} \quad \left| \sum_{j=0}^{k-1} \omega^{M}_{j,k} (b_j + \beta F_j) a_p^j \right|  \leq P
$$
hence
$$
\left| \sum_{j=0}^{k-1} \omega^{M}_{j,k} F_j a_p^j \right|  \leq \frac{2P}{\alpha - \beta}.
$$
This contradicts the property (\ref{thm:Baire:eq30}) of $\mathbf{F}$. The conclusion follows.
\qed\enddemo

Theorem \ref{thm:Baire} and Theorem \ref{thm:prevalence} mean that the image $j^{\infty}(\mathcal{E}^{0,r}_{\{M\}}) \cap \Lambda^r_{(N)}$ is generically small in the space $\Lambda^r_{(N)}$ for {\itshape any} given quasianalytic sequences $M$ and $N$ such that $\lim_{k\rightarrow + \infty}(n_k)^{1/k}=+\infty$.\vspace{6pt}

In the Roumieu case $\Lambda^r_{\{N\}}$, the notions of genericity previously used are not well defined. One can however wonder if the image is also ``small'' and in what sense. Following Remark \ref{rem:Baire}, a first direction is to obtain that the complement of the image is dense. A second possibility is to use the notion of {\itshape lineability}.

\begin{theorem}\label{thm:dense}
Let $M$ and $N$ be two quasianalytic weight sequences.\\
Let us assume that $\sup_{k\in\NN_{>0}}(n_k)^{1/k}=+\infty$, i.e. $\mathcal{O}^{0,r}\subsetneq\mathcal{E}_{\{N\}}^{0,r}$. Then $\Lambda^r_{\{N\}} \setminus j^{\infty}(\mathcal{E}^{0,r}_{\{M\}})$ is dense in $\Lambda^r_{\{N\}}$ (and so $\Lambda^r_{\{N\}} \setminus j^{\infty}(\mathcal{E}^{0,r}_{(M)})$ too).
\end{theorem}

\demo{Proof}
We proceed as in the proof of Theorem \ref{thm:Baire}. Let $\mathbf{b} \in \Lambda^r_{\{N\}}$ and let $\mathbf{F}  \in \Lambda^r_{\{N\}}$ denote a sequence whose restriction $F_j = F_{(j,0 \dots, 0)}$, $j \in \NN$, is given by Lemma \ref{lem_Thilliez}. For any $\varepsilon >0$, either $\mathbf{b} + \varepsilon \mathbf{F}$ or $\mathbf{b} - \varepsilon \mathbf{F}$ belongs to $\Lambda^r_{\{N\}} \setminus j^{\infty}(\mathcal{E}^{0,1}_{\{M\}})$. Indeed, otherwise $2\varepsilon \mathbf{F} $ would belong to $ j^{\infty}(\mathcal{E}^{0,r}_{\{M\}})$, hence also $\mathbf{F}$ which is impossible by Corollary~\ref{cor_Thilliez}. Let $h_1,h_2>0$ be such that
$$
|\mathbf{F}|^{N}_{h_1} < + \infty \text{ and } |\mathbf{b}|^{N}_{h_2}< + \infty .
$$
If $h = \max \{ h_1, h_2\}$, one has $|\mathbf{F}|^{N}_{h} < + \infty$ and $|\mathbf{b}|^{N}_{h} < + \infty$, and
$$
|\mathbf{b}-(\mathbf{b} \pm \varepsilon \mathbf{F}) |^{N}_{h} \leq \varepsilon |\mathbf{F}|^{N}_{h}.
$$
The conclusion follows.
\qed\enddemo

Let us now concentrate on the notion of {\itshape lineability}.

\begin{theorem}\label{thm:lineab}
Let $M$ and $N$ be two quasianalytic weight sequences.

Let us assume that $\sup_{k\in\NN_{>0}}(n_k)^{1/k}=+\infty$ resp. $\lim_{k\rightarrow\infty}(n_k)^{1/k}=+\infty$ , i.e. $\mathcal{O}^{0,r}\subsetneq\mathcal{E}_{\{N\}}^{0,r}$ resp. $\mathcal{O}^{0,r}\subsetneq\mathcal{E}_{(N)}^{0,r}$.

Then $\Lambda^r_{[N]} \setminus j^{\infty}(\mathcal{E}^{0,r}_{\{M\}})$ is lineable in $\Lambda^r_{[N]}$ (and so $\Lambda^r_{[N]} \setminus j^{\infty}(\mathcal{E}^{0,r}_{(M)})$ too).
\end{theorem}

\demo{Proof}
Let $\mathbf{F} \in \Lambda^r_{[N]}$ denote a sequence whose restriction $F_j = F_{(j,0 \dots, 0)}$, $j \in \NN$, is defined using  \eqref{squenceF}. For any $\lambda>0$, we define the sequence $\mathbf{F}^\lambda$ by setting
  \begin{equation*}
F^\lambda_0 = 0 \quad \text{and} \quad F^\lambda_j = \frac{F_j}{j^\lambda} , \, \forall j \in \NN_{>0}.
 \end{equation*}
Let us show that $\mathbf{F}^\lambda \notin j^\infty(\mathcal{E}^{0,r}_{\{M\}})$. From Corollary \ref{cor_Thilliez}, it suffices to prove that
 $$
 \lim_{p \to + \infty} \left| \sum_{j=0}^{k_p-1} \omega^{M}_{j,k_p}F^{\lambda}_j a^j \right|= + \infty
 $$
 for every $a$ small enough. Remark that
$$
\lim_{p \to + \infty}  \sum_{j=0}^{k_p-1} \omega^{M}_{j,k_p} F^{\lambda}_j a^j = \lim_{p\rightarrow + \infty} \left(\sum_{q=0}^{p-1}\frac{F_{k_q}}{k_q^\lambda}a^{k_q} + \sum_{j=1}^{k_{p-1}}(\omega^{M}_{j,k_p}-1)\frac{F_j}{j^\lambda}a^j \right)
 $$
 for all $a$ small enough. By (\ref{sequencelequ2}), the second term of the sum is bounded uniformly by $1$, while the first one is divergent using (\ref{powerseriesrad}).

Let $\mathcal{S}$ denote the subspace of $\Lambda^r_{[N]}$ spanned by the  $\mathbf{F}^\lambda$, $\lambda >0$: the elements of $\mathcal{S}$ can be written as
  $$
\sum_{l=1}^L \alpha_l \mathbf{F}^{\lambda_l}
 $$
 for some $L \in \NN_{>0}$, $\lambda_1, \dots, \lambda_L \in \RR$ and $\alpha_1, \dots, \alpha_L \in \CC$. Let us prove that $\mathcal{S}\setminus\{0\} \subset \Lambda _{[N]} \setminus j^{\infty}(\mathcal{E}^{0,r}_{\{M\}})$. So, let us consider $L \in \NN_{>0}$, $0<\lambda_1< \dots < \lambda_L$ in $\RR$, $\alpha_1, \dots, \alpha_L \in \CC\setminus \{0\}$, and
 $$
 \mathbf{G} = \sum_{l=1}^L \alpha_l \mathbf{F}^{\lambda_l} .
 $$
For any $a\in (0,1)$ small enough, we have
 \begin{equation}\label{eq_lineab1}
\limsup_{k\rightarrow+ \infty}\sum_{j=0}^{k-1}\omega^{M}_{j,k} G_j a^j \geq \limsup_{p \to + \infty} \sum_{j=0}^{k_p-1}\omega^{M}_{j,k_p} G_j a^j = \limsup_{p \to +\infty} \sum_{j=0}^{k_{p-1}}\omega^{M}_{j,k_p} G_j a^j
\end{equation}
 since $G_j = 0$ if $j \notin\{k_p: p\in\NN\}$. Note that
 \begin{equation}\label{eq_lineab2}
 \sum_{j=0}^{k_{p-1}}\omega^{M}_{j,k_p} G_j a^j = \sum_{j=0}^{k_{p-1}}(\omega^{M}_{j,k_p}- 1) G_j a^j + \sum_{q=0}^{p-1} G_{k_q} a^{k_q}
 \end{equation}
 and the first term of this sum can be bounded as follows uniformly in $p$
 $$
 \left| \sum_{j=0}^{k_{p-1}}(\omega^{M}_{j,k_p}- 1) G_j a^j \right| \leq \sum_{j=0}^{k_{p-1}}|\omega^{M}_{j,k_p}- 1| |F_j|\big( |\alpha_1| + \dots + |\alpha_L| \big)  \leq \big( |\alpha_1| + \dots + |\alpha_L| \big),
 $$
 by using (\ref{sequencelequ2}). However, the partial sums of the power series
 $$
 \sum_{q=0}^{+ \infty} G_{k_q} a^{k_q}
 $$
cannot be bounded: indeed, otherwise one would have
 $$
 \limsup_{q \to + \infty} \left|G_{k_q}\right|^{\frac{1}{k_q}} < + \infty \, ,
 $$
and noting that
 \begin{equation*}\label{eq:lineability_lim}
  \lim_{j \to + \infty} \frac{G_j j^{\lambda_1}}{F_j} = \lim_{j \to + \infty} \alpha_1+\alpha_2j^{\lambda_1-\lambda_2}+\dots+\alpha_l j^{\lambda_1-\lambda_L} = \alpha_1,
 \end{equation*}
 this would in turn imply that
  $$
 \limsup_{q \to + \infty} \left|F_{k_q}\right|^{\frac{1}{k_q}} < + \infty \, ,
 $$
 which is impossible from the choice of $\mathbf{F}$ (see \eqref{squenceF} and \eqref{powerseriesrad}). Hence, using \eqref{eq_lineab2}, one has
 $$
 \limsup_{p \to + \infty} \left| \sum_{j=0}^{k_{p-1}}\omega^{M}_{j,k_p} G_j a^j \right|= + \infty,$$
 and together with \eqref{eq_lineab1}, Corollary \ref{cor_Thilliez} and Remark \ref{rem:dim_r}, it gives that  $\mathbf{G}$ does not belong to the image of the Borel map~$j^\infty(\mathcal{E}^{0,r}_{\{M\}})$. The conclusion follows.
 \qed\enddemo

\begin{remark}
The vector subspace constructed in Theorem \ref{thm:lineab} has a maximal dimension which is the dimension  of the set  $\{ \lambda : \lambda >0 \}$ in $\Lambda^r_{[N]}$. We say that $\Lambda^r_{[N]} \setminus j^{\infty}(\mathcal{E}^{0,r}_{\{M\}})$ is maximal-lineable in $\Lambda^r_{[N]}$.
\end{remark}

\section{The weight matrix case $\mathcal{M}$}\label{matrixcase}
\subsection{General definitions}

\begin{definition}
A {\itshape weight matrix} $\mathcal{M}$ associated with $\RR_{>0}$ is a (one parameter) family of sequences $\mathcal{M}:=\{M^{(\lambda)}\in\RR_{>0}^{\NN}: \lambda\in\RR_{>0}\}$, such that
$$\forall\;\lambda\in\RR_{>0},\;\;M^{(\lambda)}\;\text{is a weight sequence}$$
and
$$\;M^{(\lambda)}\le M^{(\kappa)} \;\; \text{(which is equivalent to } m^{(\lambda)}\le m^{(\kappa)})\;\text{for all}\;\lambda\le\kappa,
$$
where we have put $m^{(\lambda)}_p:=\frac{M^{(\lambda)}_p}{p!}$ for $p\in\NN$.

A matrix is called {\itshape constant} if $M^{(\lambda)}\hyperlink{approx}{\approx}M^{(\kappa)}$ (i.e. $M^{(\lambda)} \hyperlink{mpreceq}{\preceq} M^{(\kappa)}$ and $M^{(\lambda)} \hyperlink{mpreceq}{\preceq} M^{(\kappa)}$) for all $\lambda,\kappa\in\RR_{>0}$.\vspace{6pt}
\end{definition}

We introduce classes of ultradifferentiable function of Roumieu type $\mathcal{E}_{\{\mathcal{M}\}}$ and of Beurling type $\mathcal{E}_{(\mathcal{M})}$ as follows, see \cite[Section 7]{dissertation} and \cite[Section 4.2]{compositionpaper}.

\begin{definition}
Let $r\in\NN_{>0}$ and $U\subseteq\RR^r$ be non-empty and open. The \emph{$\mathcal{M}$-ultradifferentiable classes of Roumieu and Beurling types} are defined respectively by
\begin{equation*}\label{generalroumieu}
\mathcal{E}_{\{\mathcal{M}\}}(U):=\bigcap_{K\subseteq U}\bigcup_{\lambda\in\RR_{>0}}\mathcal{E}_{\{M^{(\lambda)}\}}(K)
\end{equation*}
and
\begin{equation*}\label{generalbeurling}
\mathcal{E}_{(\mathcal{M})}(U):=\bigcap_{\lambda\in\RR_{>0}}\mathcal{E}_{(M^{(\lambda)})}(U).
\end{equation*}
\end{definition}

For a compact set $K\subseteq\RR^r$, one has the representations
$$\mathcal{E}_{\{\mathcal{M}\}}(K):=\underset{\lambda\in\RR_{>0}}{\varinjlim}\;\underset{h>0}{\varinjlim}\;\mathcal{E}_{M^{(\lambda)},h}(K)$$
and so for $U\subseteq\RR^r$ non-empty open
\begin{equation*}\label{generalroumieu1}
\mathcal{E}_{\{\mathcal{M}\}}(U)=\underset{K\subseteq U}{\varprojlim}\;\underset{\lambda\in\RR_{>0}}{\varinjlim}\;\underset{h>0}{\varinjlim}\;\mathcal{E}_{M^{(\lambda)},h}(K).
\end{equation*}
Similarly we get for the Beurling case
\begin{equation*}\label{generalbeurling1}
\mathcal{E}_{(\mathcal{M})}(U)=\underset{K\subseteq U}{\varprojlim}\;\underset{\lambda\in\RR_{>0}}{\varprojlim}\;\underset{h>0}{\varprojlim}\;\mathcal{E}_{M^{(\lambda)},h}(K).
\end{equation*}
Consequently, since the sequences of $\mathcal{M}$ are pointwise ordered, $\mathcal{E}_{(\mathcal{M})}(U)$ is a {\itshape Fr\'echet space} and

$\underset{\lambda\in\RR_{>0}}{\varinjlim}\;\underset{h>0}{\varinjlim}\;\mathcal{E}_{M^{(\lambda)},h}(K)=\underset{n\in\NN_{>0}}{\varinjlim}\;\;\mathcal{E}_{M^{(n)},n}(K)$ is a {\itshape Silva space}, i.e. a countable inductive limit of Banach spaces with compact connecting mappings. For more details concerning the locally convex topology in this setting we refer to \cite[Section 4.2]{compositionpaper}.\vspace{6pt}

\begin{definition}
The spaces of \emph{germs at $0\in\RR^r$ of the $(\mathcal{M})$-ultradifferentiable functions of Roumieu and Beurling types} are defined respectively by
\begin{equation*}\label{roumieumatrixgerm}
\mathcal{E}_{\{\mathcal{M}\}}^{0,r}:=\underset{k\in\NN_{>0}}{\varinjlim}\mathcal{E}_{\{\mathcal{M}\}}\left(\left(-\frac{1}{k},\frac{1}{k}\right)^r\right),
\end{equation*}
and
\begin{equation*}\label{beurlingmatrixgerm}
\mathcal{E}_{(\mathcal{M})}^{0,r}:=\underset{k\in\NN_{>0}}{\varinjlim}\mathcal{E}_{(\mathcal{M})}\left(\left(-\frac{1}{k},\frac{1}{k}\right)^r\right).
\end{equation*}
\end{definition}

Finally, as done in the case of weight sequences, we introduce the corresponding spaces of sequences, and we endow them with their classical topology.

\begin{definition}
We introduce the sequence classes of Roumieu type
$$\Lambda^r_{\{\mathcal{M}\}}:=\bigcup_{\lambda\in\RR_{>0}}\Lambda^r_{\{ M^{(\lambda)}\} }=\left\{\mathbf{b}\in\CC^{\NN^r} : \exists\;\lambda\in\RR_{>0}\;\;\exists h > 0 , |\mathbf{b}|^{M^{(\lambda)}}_{h}< + \infty \right\},$$
and of Beurling type
$$\Lambda^r_{(\mathcal{M})}:=\bigcap_{\lambda\in\RR_{>0}}\Lambda^r_{(M^{(\lambda)}) }=\left\{\mathbf{b}\in\CC^{\NN^r} : \forall\;\lambda\in\RR_{>0}\;\;\forall h > 0 , |\mathbf{b}|^{M^{(\lambda)}}_{h}< + \infty \right\}.$$
\end{definition}

Using notations similar as before, the \emph{Borel map} $j^\infty$ is defined in the weight matrix case by
\begin{equation*}\label{Borelmapmatrix}
j^{\infty}:\mathcal{E}_{[\mathcal{M}]}^{0,r}\longrightarrow\Lambda^r_{[\mathcal{M}]},\hspace{20pt}j^\infty(f) = \left(\frac{\partial^{\alpha}f(0)}{|\alpha|!}\right)_{\alpha\in \NN^r}.
\end{equation*}

In \cite[Theorem 4.1]{testfunctioncharacterization}, the following result has been obtained (under slightly more general assumptions on $\mathcal{M}$ and using regularizations of $M^{(\lambda)}$).

\begin{theorem}\label{Matrix-non-quasi-analyticity}
Let $\mathcal{M}=\{M^{(\lambda)}: \lambda\in\RR_{>0}\}$ be a weight matrix.
\begin{itemize}
\item[$(i)$] $\mathcal{E}_{\{\mathcal{M}\}}$ is non-quasianalytic if and only if there exists $\lambda_0\in\RR_{>0}$ such that $\mathcal{E}_{[M^{(\lambda_0)}]}$ is non-quasianalytic.

\item[$(ii)$] $\mathcal{E}_{(\mathcal{M})}$ is non-quasianalytic if and only if each $\mathcal{E}_{[M^{(\lambda)}]}$ is non-quasianalytic.
\end{itemize}
\end{theorem}

This result yields and motivates the following definition, see also \cite[Section 5.1]{borelmappingquasianalytic}.

\begin{definition}\label{matrixquasianalytic}
A weight matrix $\mathcal{M}$ is called \emph{quasianalytic} if for all $\lambda\in\RR_{>0}$ the sequence $M^{(\lambda)}$ is quasianalytic, which means
\begin{equation}\label{matrixquasianalyticequ}
\forall\;\lambda\in\RR_{>0},\;\;\;\sum_{j=1}^{+\infty}\frac{1}{\left(M^{(\lambda)}_j\right)^{1/j}}=+\infty.
\end{equation}
\end{definition}

In this case both classes $\mathcal{E}_{\{\mathcal{M}\}}$ and $\mathcal{E}_{(\mathcal{M})}$ and all classes $\mathcal{E}_{\{M^{(\lambda)}\}}$ resp. $\mathcal{E}_{(M^{(\lambda)})}$ are quasianalytic too, see Proposition \ref{nonquasiremarks}. For the Beurling case $\mathcal{E}_{(\mathcal{M})}$ it would be enough to require only that there is some $M^{(\lambda_0)}$ which is quasianalytic since then $M^{(\lambda)}$ for all $\lambda\le\lambda_0$ is quasianalytic too and since, by definition of the Beurling type classes, the spaces remain unchanged if we remove from $\mathcal{M}$ all (possible non-quasianalytic sequences) $M^{(\lambda)}$ for $\lambda>\lambda_0$.\vspace{6pt}

Recently, in \cite[Thm. 5, Thm. 6]{borelmappingquasianalytic}, it has been shown that $j^{\infty}$ restricted to the germs $\mathcal{E}_{[\mathcal{M}]}^{0,r}$ can never be onto the corresponding sequence space for any quasianalytic weight matrix $\mathcal{M}$ such that $\mathcal{O}^{0,r}\subsetneq\mathcal{E}_{[\mathcal{M}]}^{0,r}$.

\subsection{Generalization of Thilliez's proof for non-surjectivity}\label{generalizationThilliezmatrix}
We generalize \cite[Theorem 3]{thilliez} resp. Theorem \ref{Thillieztheorem} to the general weight matrix setting (for both types). To ensure $\mathcal{O}^{0,r}\subsetneq\mathcal{E}_{\{\mathcal{M}\}}^{0,r}$ resp. $\mathcal{O}^{0,r}\subsetneq\mathcal{E}_{(\mathcal{M})}^{0,r}$ we assume
$$\forall\;\lambda\in\RR_{>0},\;\;\;\sup_{k\in\NN_{>0}}\left(m^{(\lambda)}_k\right)^{1/k}=+\infty\hspace{25pt}\text{resp.}\hspace{25pt}\forall\;\lambda\in\RR_{>0},\;\;\;\lim_{k\rightarrow+ \infty}\left(m^{(\lambda)}_k\right)^{1/k}=+\infty,$$ e.g. see \cite[Section 5]{borelmappingquasianalytic} and which follows from \cite[Proposition 4.6]{compositionpaper}. Note that in the Roumieu case, one could assume that $\sup_{k\in\NN_{>0}}\left(m^{(\lambda)}_k\right)^{1/k}=+\infty$ only for all $\lambda\ge\lambda_0$ for some $\lambda_0\in\RR_{>0}$ (large): Indeed, one can skip in this case all small sequences in the matrix without changing the ultradifferentiable class.

Let us start with a generalization of Lemma \ref{lem_Thilliez} working with sequence spaces defined via weight matrices.

\begin{lemma}\label{lem_Thilliezmatrix}
Let $M$ be a quasianalytic weight sequence and $\mathcal{N}=\{N^{(\lambda)}: \lambda\in\RR_{>0}\}$ be a quasianalytic weight matrix such that $\sup_{k\in\NN_{>0}}\left(n^{(\lambda)}_k\right)^{1/k}=+\infty$ resp. $\lim_{k\rightarrow + \infty}\left(n^{(\lambda)}_k\right)^{1/k}=+\infty$ for all $\lambda\in\RR_{>0}$, i.e. $\mathcal{O}^{0,1}\subsetneq\mathcal{E}_{\{\mathcal{N}\}}^{0,1}$ resp. $\mathcal{O}^{0,1}\subsetneq\mathcal{E}_{(\mathcal{N})}^{0,1}$. There exists $\mathbf{F} \in \Lambda^1_{[\mathcal{N}]}$ and $a_0  \in (0,1]$ such that
\begin{equation*}\label{eq_cor_Thilliez2}
\limsup_{k \to + \infty} \left| \sum_{j=0}^{k-1} \omega^{M}_{j,k} F_j a^j \right|= + \infty,
\end{equation*}
for all $0<a \leq a_0$.
\end{lemma}

\demo{Proof}
Note first that the Roumieu case follows immediately from Lemma \ref{lem_Thilliez}: indeed, it suffices to fix $N^{(\lambda_0)}\in \mathcal{N}$ and to use the inclusion $\Lambda^1_{\{N^{(\lambda_0)}\}}\subseteq\Lambda^1_{\{\mathcal{N}\}}$.
So, let us concentrate on the Beurling case. Let $0<a_0\le 1$ be arbitrary but from now on fixed. Let us show that there is a strictly increasing sequence $(k_p)_{p\in \NN}$ with  $k_0\ge 1$ such that the sequence $\mathbf{F}$ defined by
\begin{equation}\label{squenceFbeur}
F_{k_p}:=\sqrt{n^{\big(\frac{1}{p+1}\big)}_{k_p}} \quad \text{and} \quad F_j :=  0 \, \, \text{ otherwise, }
\end{equation}
satisfies
\begin{equation}\label{sequencelequ2matrix}
\sum_{j=0}^{k_{p-1}}\left|\omega^M_{j,k_p}-1\right||F_j|a_0^j\le 1, \quad \forall p \in \NN_{>0},
\end{equation}
and
\begin{equation}\label{powerseriesradmatrix}
\lim_{p\rightarrow + \infty}F_{k_p}^\frac{1}{k_p}=+\infty.
\end{equation}
It is easy to show this by induction, since $\lim_{k \to + \infty} \omega^M_{j,k}=1$ for all $j \in \NN$ by \eqref{eq:lim_omega}. So if $k_{p-1}$ has been constructed, one may choose $k_p$ sufficiently large in order to have both \eqref{sequencelequ2matrix} and
$$
F_{k_p}^\frac{1}{k_p} = \sqrt{\left(n_{k_p}^{\big(\frac{1}{p+1}\big)}\right)^{\frac{1}{k_p}}} \geq p.
$$
This is possible by using the assumption $\lim_{k\rightarrow + \infty}\left(n^{(\lambda)}_k\right)^{1/k}=+\infty$ for all $\lambda\in\RR_{>0}$ and guarantees \eqref{powerseriesradmatrix}.

That $\mathbf{F}$ belongs to $\Lambda_{(\mathcal{N})}^1$ follows easily from the definition and the assumption of the lemma. Moreover, for all $a \in (0,a_0]$, one has
$$
\sum_{j=0}^{k_p-1} \omega^M_{j,k_p} F_j a^j = \sum_{q=0}^{p-1} F_{k_q} a^{k_q} + \sum_{j=0}^{k_p-1} \big(\omega^M_{j,k_p}-1\big) F_j a^j
$$
which implies
$$
\limsup_{k \to + \infty} \sum_{j=0}^{k-1} \omega^M_{j,k} F_j a^j = \infty
$$
by \eqref{sequencelequ2matrix} and \eqref{powerseriesradmatrix}.
 \qed\enddemo

\begin{remark} As defined in \eqref{squenceFbeur}, the sequence $\mathbf{F}$ can never define a real analytic germ: Otherwise, $F_{k_p}=\sqrt{n^{(1/(p+1))}_{k_p}}\le C h^{k_p}$ should be satisfied for some $C,h>0$ and all $p\in\NN$, a contradiction to \eqref{powerseriesradmatrix}. The Roumieu case follows immediately by Lemma \ref{lem_Thilliez} and \eqref{powerseriesrad} for $N^{(\lambda_0)}$. Moreover, as commented in Remark \ref{remarklpj}, in the proof of Lemma \ref{lem_Thilliezmatrix} the sequence $\mathbf{F}$ could be defined by setting $F_j := l^p_j$ for $k_{p-1}<j<k_p$ for non-zero values $l^p_j$ subjected to some precise growth control. This provides some additional information on sequences not contained in the image of the Borel map but is not necessary for the forthcoming proofs.
\end{remark}

By proceeding as in the weight sequence case, Corollary \ref{cor_Thilliez} directly gives the following proposition.

\begin{proposition}\label{prop_matrix}
Let $M$ be a quasianalytic weight sequence and $\mathcal{N}=\{N^{(\lambda)}: \lambda\in\RR_{>0}\}$ be a quasianalytic weight matrix such that $\sup_{k\in\NN_{>0}}\left(n^{(\lambda)}_k\right)^{1/k}=+\infty$ resp. $\lim_{k\rightarrow + \infty}\left(n^{(\lambda)}_k\right)^{1/k}=+\infty$ for all $\lambda\in\RR_{>0}$, i.e. $\mathcal{O}^{0,r}\subsetneq\mathcal{E}_{\{\mathcal{N}\}}^{0,r}$ resp. $\mathcal{O}^{0,r}\subsetneq\mathcal{E}_{(\mathcal{N})}^{0,r}$. Then, one has $$j^{\infty}(\mathcal{E}_{\{M\}}^{0,r}) \cap \Lambda^r_{[\mathcal{N}]} \subsetneq \Lambda^r_{[\mathcal{N}]}$$ (and hence $j^{\infty}(\mathcal{E}_{(M)}^{0,r}) \cap \Lambda^r_{[\mathcal{N}]} \subsetneq \Lambda^r_{[\mathcal{N}]}$ also).
\end{proposition}

Note that in order to get the non-surjectivity of the Borel map in the weight matrix case, we need to get an equivalent of Proposition \ref{prop_matrix} working only with quasianalytic weight matrices (and not with a weight sequence). This can be obtained thanks to the following result.

\begin{proposition}\label{Thillieztheoremmatrixgeneralizationroum}
Let $\mathcal{M}=\{M^{(\lambda)}: \lambda\in\RR_{>0}\}$ be a quasianalytic weight matrix. Then there exists a quasianalytic weight sequence $L$ satisfying $M^{(\lambda)}\hyperlink{mtriangle}{\vartriangleleft}L$ for all $\lambda\in\RR_{>0}$, i.e. $\mathcal{E}_{\{\mathcal{M}\}} \subseteq \mathcal{E}_{(L)}$ holds true.
\end{proposition}

The aim is to construct a quasianalytic (weight) sequence $L$ lying (strictly) above $\mathcal{M}$ by applying some diagonal technique. Unfortunately, it seems that such a construction does not preserve the log-convexity; we can overcome this problem by working with regularizations of $L$ and by applying Proposition \ref{nonquasiremarks}. The following idea is motivated by the proof of \cite[Prop. 4.7 (i)]{testfunctioncharacterization}.\vspace{6pt}

\demo{Proof}
Let $(d_i)_{i\in \NN_{>0}}$ be a strictly increasing sequence in $\RR$, with $d_1\ge 1$ and tending to infinity as $i\rightarrow+\infty$. By the assumptions on $\mathcal{M}$ there exists a strictly increasing sequence $(j_i)_{i\in \NN_{>0}}$ (in $\NN$) with $j_1=1$ and such that $\sum_{j=j_i}^{j_{i+1}-1}\frac{1}{(M^{(i)}_j)^{1/j}}\ge d_i$, see \eqref{matrixquasianalyticequ}. According to this sequence, we put
$$\widetilde{L}_0=1 \quad \text{and} \quad\widetilde{L}_j =  d^j_i M^{(i)}_j, \quad \text{if}\;j_i\le j<j_{i+1},\;\forall i \in \NN_{>0}.$$
First, for any given index $\lambda_0\in\RR_{>0}$ (large), we have $M^{(\lambda_0)}_j\leq M^{(i)}_j=d_i^{-j}\widetilde{L}_j$ for all $j_i\le j<j_{i+1}$ and $i\ge\lambda_0$, hence $M^{(\lambda)}\hyperlink{mtriangle}{\vartriangleleft}\widetilde{L}$ for all $\lambda\in\RR_{>0}$ because $d_i\rightarrow+ \infty$ as $i\rightarrow+ \infty$. So $\mathcal{E}_{\{M^{(\lambda_0)}\}} \subseteq \mathcal{E}_{(\widetilde{L})}$ and $\mathcal{E}_{\{\mathcal{M}\}} \subseteq \mathcal{E}_{(\widetilde{L})}$ follows by definition as a special case by \cite[Prop. 4.6 (2)]{compositionpaper}.

Unfortunately we do not see directly if $\widetilde{L}$ is log convex but since $d_i \rightarrow + \infty$ as $i \to + \infty$, and since $0<\liminf_{p\rightarrow+ \infty}\left(m^{(\lambda)}_p\right)^{1/p}$ for all $\lambda\in\RR_{>0}$, we obtain $\lim_{p\rightarrow+ \infty}(\widetilde{l}_p)^{1/p}=+\infty$ (where $\widetilde{l}_p=\frac{\widetilde{L}_p}{p!}$). Consequently $\mathcal{E}_{[\widetilde{L}]}=\mathcal{E}_{[\widetilde{L}^{\on{lc}}]}$, i.e. $\widetilde{L}$ can be replaced by its log-convex minorant for both cases, see Remark \ref{logminorantremark}, and so $M^{(\lambda)}\hyperlink{mtriangle}{\vartriangleleft}\widetilde{L}^{\on{lc}}$ for all $\lambda \in\RR_{>0}$ too.

It remains to show that (\ref{eq_Q}) holds true for $\widetilde{L}^{\on{lc}}$. By definition of $\widetilde{L}$, the log-convexity of each $M^i$ and since $(d_i)_i$ is (strictly) increasing, we have that $j\mapsto(\widetilde{L}_j)^{1/j}$ is increasing: If $j=j_{i+1}$, then for all $i\ge 1$, one has $$(\widetilde{L}_{j_{i+1}-1})^{1/(j_{i+1}-1)}=d_i(M^{(i)}_{j_{i+1}-1})^{1/(j_{i+1}-1)}\le d_i(M^{(i)}_{j_{i+1}})^{1/j_{i+1}}\le d_{i+1}(M^{(i+1)}_{j_{i+1}})^{1/j_{i+1}}=(\widetilde{L}_{j_{i+1}})^{1/j_{i+1}}.$$
The remaining cases are clear.

So we have shown $\widetilde{L}=\widetilde{L}^I$ and finally $$\sum_{j=1}^{+\infty}\frac{1}{(\widetilde{L}^I_j)^{1/j}}=\sum_{i=1}^{+\infty}\sum_{j=j_i}^{j_{i+1}-1}\frac{1}{(\widetilde{L}_j)^{1/j}}=\sum_{i=1}^{+ \infty}\sum_{j=j_i}^{j_{i+1}-1}\frac{1}{d_i(M^{(i)}_j)^{1/j}}\ge\sum_{i=1}^{+\infty}1=+\infty,$$ by the choice of $(j_i)_{i \in \NN_{>0}}$ above. Hence by Proposition \ref{nonquasiremarks} we get (\ref{eq_Q}) for $\widetilde{L}^{\on{lc}}$ and the conclusion follows by taking $L:=\widetilde{L}^{\on{lc}}$.
\qed\enddemo

\medskip

Propositions \ref{prop_matrix} and  \ref{Thillieztheoremmatrixgeneralizationroum} imply directly the next main result.

\begin{theorem}\label{Thillieztheoremmatrixgeneralization}
Let $\mathcal{M}$ and $\mathcal{N}$ be two quasianalytic weight matrices such that $\sup_{k\in\NN_{>0}}\left(n^{(\lambda)}_k\right)^{1/k}=+\infty$ resp. $\lim_{k\rightarrow+ \infty}\left(n^{(\lambda)}_k\right)^{1/k}=+\infty$ for all $\lambda\in\RR_{>0}$, i.e. $\mathcal{O}^{0,r}\subsetneq\mathcal{E}_{\{\mathcal{N}\}}^{0,r}$ resp. $\mathcal{O}^{0,r}\subsetneq\mathcal{E}_{(\mathcal{N})}^{0,r}$. Then, one has $$j^{\infty}(\mathcal{E}_{\{\mathcal{M}\}}^{0,r}) \cap \Lambda^r_{[\mathcal{N}]} \subsetneq \Lambda^r_{[\mathcal{N}]}$$ (and hence $j^{\infty}(\mathcal{E}_{(\mathcal{M})}^{0,r}) \cap \Lambda^r_{[\mathcal{N}]} \subsetneq \Lambda^r_{[\mathcal{N}]}$ also). In particular, the Borel map $j^{\infty}: \mathcal{E}_{[\mathcal{N}]}^{0,r}\longrightarrow\Lambda^r_{[\mathcal{N}]}$ is not surjective.
\end{theorem}

\subsection{Generic size of the image of the Borel map}\label{Borelmapfarawaymatrix}

Thanks to Theorem \ref{Thillieztheoremmatrixgeneralization}, we can immediately transfer and generalize all central statements from Section \ref{Borelmapfaraway} to the weight matrix case.

\begin{theorem}\label{thm:Bairematrix}
Let $\mathcal{N}$ and $\mathcal{M}$ be two quasianalytic weight matrices. Assume that $\lim_{k\rightarrow+ \infty}\left(n^{(\lambda)}_k\right)^{1/k}=+\infty$ for all $\lambda\in\RR_{>0}$, i.e. $\mathcal{O}^{0,r}\subsetneq\mathcal{E}_{(\mathcal{N})}^{0,r}$. Then the image of the Borel map $j^{\infty}(\mathcal{E}^{0,r}_{\{\mathcal{M}\}}) \cap \Lambda^r_{(\mathcal{N})}$ is meager in $\Lambda^r_{(\mathcal{N})}$.
\end{theorem}

\demo{Proof}
By Proposition \ref{Thillieztheoremmatrixgeneralizationroum}, we know that there is a quasianalytic weight sequence $L$ such that $\mathcal{E}^{0,r}_{\{\mathcal{M}\}} \subset \mathcal{E}^{0,r}_{(L)} \subset \mathcal{E}^{0,r}_{\{L\}}$. So, it suffices to prove that the larger set $j^{\infty}(\mathcal{E}^{0,r}_{\{L\}}) \cap \Lambda^r_{(\mathcal{N})}$ is meager in $\Lambda^r_{(\mathcal{N})}$. By Lemma \ref{lem_Thilliezmatrix}, we can consider $\mathbf{F} \in \Lambda^r_{(\mathcal{N})}$ and $a_0  \in (0,1]$ such that
\begin{equation*}\label{eq_cor_Thilliez2}
\limsup_{k \to + \infty} \left| \sum_{j=0}^{k-1} \omega^{L}_{j,k} F_j a^j \right|= + \infty,
\end{equation*}
for all $0<a \leq a_0$, where  $F_j = F_{(j,0,\dots,0)}$ for any $j \in \NN$. As done in the proof of Theorem \ref{thm:Baire} and using Corollary \ref{cor_Thilliez}, if $ (a_p)_{p \in \NN}$ is a fixed sequence of $(0,a_0]$ which  decreases to $0$, it suffices to prove that the set
$$
\mathcal{G} = \bigcap_{p \in \NN} \bigcap_{P \in \NN} \bigcap_{K \in \NN_0} \bigcup_{k \ge K} \left\{\mathbf{b} \in \Lambda^r_{(\mathcal{N})} : \left| \sum_{j=0}^{k-1} \omega^{L}_{j,k} b_j a_p^j \right| > P  \right\}
$$
is a countable intersection of dense open sets of $\Lambda^r_{(\mathcal{N})}$. Since the inclusion $\Lambda^r_{(\mathcal{N})} \hookrightarrow \Lambda^r_{(N^{(\lambda)})}$ is continuous for any $\lambda \in \RR_{>0}$, and from the proof of Theorem \ref{thm:Baire}, it is clear that $\mathcal{G}$ is a countable intersection of open sets. We obtain the density of these sets noting that the equality \eqref{thm:Baire:eq3} holds true for all $h>0$ and all $N^{(\lambda)}$, $\lambda\in\RR_{>0}$.
\qed\enddemo

\medskip

Similarly, we get the generalization of Theorem \ref{thm:prevalence} to the matrix setting.

\begin{theorem}\label{thm:prevalencematrix}
Let $\mathcal{M}$ and $\mathcal{N}$ be two quasianalytic weight matrices. Assume that $\lim_{k\rightarrow + \infty}\left(n^{(\lambda)}_k\right)^{1/k}=+\infty$ for all $\lambda\in\RR_{>0}$, i.e. $\mathcal{O}^{0,r}\subsetneq\mathcal{E}_{(\mathcal{N})}^{0,r}$. Then the image of the Borel map $j^{\infty}(\mathcal{E}^{0,r}_{\{\mathcal{M}\}}) \cap \Lambda^r_{(\mathcal{N})}$ is Haar-null in $\Lambda^r_{(\mathcal{N})}$.

\end{theorem}

\demo{Proof}
As done before, using Proposition \ref{Thillieztheoremmatrixgeneralizationroum}, we can reduce the proof to the case where the weight matrix $\mathcal{M}$ is constant. We follow then the lines of the proof of Theorem \ref{thm:prevalence}, where the set $\mathcal{G}$ is defined as in the proof of Theorem \ref{thm:Bairematrix}.
\qed\enddemo

\medskip

Moreover, in the Roumieu case, we have the following result.

\begin{theorem}\label{thm:densematrix}
Let $\mathcal{M}$ and $\mathcal{N}$ be two quasianalytic weight matrices. Assume that $\sup_{k\in\NN_{>0}}\left(n^{(\lambda)}_k\right)^{1/k}=+\infty$ for all $\lambda\in\RR_{>0}$, i.e. $\mathcal{O}^{0,r}\subsetneq\mathcal{E}_{\{\mathcal{N}\}}^{0,r}$. Then, the set $\Lambda^r_{\{\mathcal{N}\}} \setminus j^{\infty}(\mathcal{E}^{0,r}_{\{\mathcal{M}\}})$ is dense in $\Lambda^r_{\{\mathcal{N}\}}$.

\end{theorem}

\demo{Proof}
Again, using Proposition \ref{Thillieztheoremmatrixgeneralizationroum}, it suffices to consider the case where the weight matrix $\mathcal{M}$ is constant. In the proof of Theorem \ref{thm:dense}, let $h_1,h_2>0$ and $\lambda_1,\lambda_2>0$ be such that
$|\mathbf{F}|^{N^{(\lambda_1)}}_{h_1} < + \infty \text{ and } |\mathbf{b}|^{N^{(\lambda_2)}}_{h_2}< + \infty$. Then we put $\kappa:=\max\{\lambda_1,\lambda_2\}$ and again $h:= \max \{ h_1, h_2\}$ to get both $|\mathbf{F}|^{N^{(\kappa)}}_{h} < + \infty$ and $|\mathbf{b}|^{N^{(\kappa)}}_{h} < + \infty$. The conclusion follows.
\qed\enddemo

\medskip

Finally, we can obtain an equivalent of these results in the context of lineability.

\begin{theorem}\label{thm:lineabmatrix}
Let $\mathcal{M}$ and $\mathcal{N}$ be two quasianalytic weight matrices. Assume that $\sup_{k\in\NN_{>0}}\left(n^{(\lambda)}_k\right)^{1/k}=+\infty$ for all $\lambda\in\RR_{>0}$ resp. $\lim_{k\rightarrow + \infty}\left(n^{(\lambda)}_k\right)^{1/k}=+\infty$ for all $\lambda\in\RR_{>0}$, i.e. $\mathcal{O}^{0,r}\subsetneq\mathcal{E}_{\{\mathcal{N}\}}^{0,r}$ resp. $\mathcal{O}^{0,r}\subsetneq\mathcal{E}_{(\mathcal{N})}^{0,r}$. Then, the set $\Lambda^r_{[\mathcal{N}]} \setminus j^{\infty}(\mathcal{E}^{0,r}_{\{\mathcal{M}\}})$ is lineable in $\Lambda^r_{[\mathcal{N}]}$ (and so $\Lambda^r_{[\mathcal{N}]} \setminus j^{\infty}(\mathcal{E}^{0,r}_{(\mathcal{M})})$ too).

\end{theorem}

\demo{Proof}
As previously, we consider the case where the weight matrix $\mathcal{M}$ is constant and we follow simply the proof of Theorem \ref{thm:lineab}.
\qed\enddemo

We close this section with the following observation.

\begin{remark}
We have used the proofs from the single weight sequence case of Section \ref{Borelmapfaraway} and transferred them to the more general weight matrix case of this Section \ref{Borelmapfarawaymatrix}. Alternatively, one could start directly with the weight matrix setting (and give the proofs from Section \ref{Borelmapfaraway} in this general approach) and then obtain the single weight sequence case as an immediate consequence for the constant matrix $\mathcal{M}=\{M\}$.
\end{remark}

\section{The weight function case}\label{weighfunctioncase}
\subsection{General definitions}\label{weightfunctionclasses}

In this last part, we will study classes of ultradifferentiable functions defined using weight functions in the sense of Braun-Meise-Taylor, see \cite{BraunMeiseTaylor90}. As we will see, this case can be reduced to the weight matrix case. First, let us start by recalling the basic definitions.

\begin{definition}
A function $\omega:[0,+\infty)\rightarrow[0,+\infty)$ is called a \emph{weight function} if
\begin{itemize}
\item[$(i)$] $\omega$ is continuous,
\item[$(ii)$] $\omega$ is increasing,
\item[$(iii)$] $\omega(t)=0$ for all $t\in[0,1]$ (normalization, w.l.o.g.),
\item[$(iv)$] $\lim_{t\rightarrow+\infty}\omega(t)=+\infty$.
\end{itemize}
In this case, we say that $\omega$ has $\hypertarget{om0}{(\omega_0)}$.
\end{definition}

Classical additional conditions can be imposed on the considered weight functions. More precisely, let us define the following conditions:
\begin{itemize}
\item[\hypertarget{om1}{$(\omega_1)$}] $\omega(2t)=O(\omega(t))$ as $t\rightarrow+\infty$,

\item[\hypertarget{om2}{$(\omega_2)$}] $\omega(t)=O(t)$ as $t\rightarrow+ \infty$,

\item[\hypertarget{om3}{$(\omega_3)$}] $\log(t)=o(\omega(t))$ as $t\rightarrow+\infty$ ($\Leftrightarrow\lim_{t\rightarrow+\infty}\frac{t}{\varphi_{\omega}(t)}=0$),

\item[\hypertarget{om4}{$(\omega_4)$}] $\varphi_{\omega}:t\mapsto\omega(e^t)$ is a convex function on $\RR$,

\item[\hypertarget{om5}{$(\omega_5)$}] $\omega(t)=o(t)$ as $t\rightarrow+\infty$.

\end{itemize}

For convenience, we define the set
$$\hypertarget{omset1}{\mathcal{W}}:=\{\omega:[0,+ \infty)\rightarrow[0, + \infty): \omega\;\text{has}\;\hyperlink{om0}{(\omega_0)}, \hyperlink{om1}{(\omega_1)}, \hyperlink{om3}{(\omega_3)},\hyperlink{om4}{(\omega_4)}\}.$$
Note that \hyperlink{om2}{$(\omega_2)$} is sometimes also considered as a general assumption on $\omega$ (e.g. see \cite[Sect. 4.1]{borelmappingquasianalytic}) and note also that \hyperlink{om5}{$(\omega_5)$} implies \hyperlink{om2}{$(\omega_2)$}.

For $\omega\in\hyperlink{omset}{\mathcal{W}}$, we define the {\itshape Legendre-Fenchel-Young-conjugate} of $\varphi_{\omega}$ by
$$\varphi^{*}_{\omega}(x):=\sup\{x y-\varphi_{\omega}(y): y\ge 0\},\;\;\;x\ge 0.$$

\begin{definition}
Let $r\in\NN_{>0}$, $U\subseteq\RR^r$ be a non-empty open set and $\omega\in\hyperlink{omset}{\mathcal{W}}$. The $\omega$-\emph{ultradifferentiable Roumieu type class} is defined by
$$\mathcal{E}_{\{\omega\}}(U):=\{f\in\mathcal{E}(U):\;\forall\;K\subseteq U\;\text{compact}\;\; \exists\;l>0,\;\|f\|^{\omega}_{K,l}<+\infty\},$$
and the $\omega$-\emph{ultradifferentiable Beuling type class} by
$$\mathcal{E}_{(\omega)}(U):=\{f\in\mathcal{E}(U):\;\forall\;K\subseteq U\;\text{compact}\;\; \forall\;l>0,\;\|f\|^{\omega}_{K,l}<+\infty\},$$
where we have put
\begin{equation*}\label{semi-norm-1}
\|f\|^{\omega}_{K,l}:=\sup_{\alpha\in\NN^r,x\in K}\frac{|\partial^{\alpha}f(x)|}{\exp(\frac{1}{l}\varphi^{*}_{\omega}(l|\alpha|))}.
\end{equation*}
\end{definition}

As done in the previous contexts, these  spaces are endowed with their natural topologies.
Analogously as in the sections above, we also consider the spaces of germs at $0$, denoted $\mathcal{E}_{\{\omega\}}^{0,r}$ and $\mathcal{E}_{(\omega)}^{0,r}$, and the associated spaces of complex sequences $\Lambda^r_{\{ \omega\} }$ and $\Lambda^r_{(\omega) }$.
Again, we endow these spaces  with their natural topology: $\Lambda_{\{\omega\} }$ is an (LB)-space and $\Lambda_{(\omega) }$ a Fr\'echet space. In this setting, the  Borel map is given by
\begin{equation*}\label{Borelmapfunction}
j^{\infty}:\mathcal{E}_{[\omega]}^{0,r}\longrightarrow\Lambda^r_{[\omega]},\hspace{20pt}j^\infty(f) = \left(\frac{\partial^{\alpha}f(0)}{|\alpha|!}\right)_{\alpha\in \NN^r}.
\end{equation*}
As pointed out in \cite[Section 4.2]{borelmappingquasianalytic}, that to ensure $\mathcal{C}^{\omega}\subsetneq\mathcal{E}_{\{\omega\}}$ resp. $\mathcal{C}^{\omega}\subsetneq\mathcal{E}_{(\omega)}$, one has to assume that $$\liminf_{t\rightarrow+ \infty}\frac{\omega(t)}{t}=0\hspace{30pt}\text{resp.}\hspace{30pt}\omega(t)=o(t)\;\text{as}\;t\rightarrow+ \infty ,\;\text{i.e.}\;\hyperlink{om5}{(\omega_5)},$$
which follows from the characterizations given in \cite[Lemm. 5.16, Cor. 5.17]{compositionpaper} and the fact that the weight $\omega(t)=t$ (up to equivalence) defines the class $\mathcal{C}^{\omega}$.

Moreover, in the present setting, the definition of quasianalyticity takes the following form.

\begin{definition}
A weight function is called \emph{quasianalytic} if it satisfies
\begin{equation*}\tag{$\omega_{\text{Q}}$}
\int_1^{+ \infty}\frac{\omega(t)}{t^2}dt= + \infty.
\end{equation*}
\end{definition}

\subsection{Generic size of the image of the Borel map}\label{Borelmapfarawayweightfunction}
Naturally, one could wonder if the results of the previous sections concerning the Borel map still hold in the context of weight functions, and if the proofs will require some different techniques and methods in this setting. But we will see that they can be obtained without any additional work! Applying the idea presented below has already been helpful in \cite{borelmappingquasianalytic} where a closely related topic has been treated.

\medskip

In \cite{dissertation} and \cite[Section 5]{compositionpaper}, a matrix $\Omega:=\{W^{(l)}=(W^{(l)}_j)_{j\in\NN}: l>0\}$ has been associated with each $\omega\in\hyperlink{omset}{\mathcal{W}}$: This matrix is defined by
$$W^{(l)}_j:=\exp\left(\frac{1}{l}\varphi^{*}_{\omega}(lj)\right) , \quad \forall j\in\NN, \, \forall l >0,$$
and  $\mathcal{E}_{[\omega]}=\mathcal{E}_{[\Omega]}$ holds as locally convex vector spaces. Moreover, the following results have been obtained:

\begin{itemize}
\item[$(i)$] Each $W^{(l)}$ satisfies the basic assumptions $(I)$ and $(II)$.

\item[$(ii)$] $\omega$ has in addition \hyperlink{om2}{$(\omega_2)$} if and only if some/each $W^{(l)}$ has $(III)$, too.
\end{itemize}

So each $W^{(l)}\in\Omega$ is a weight sequence according to the requirements from Section \ref{weightsequences}, provided $\omega\in\hyperlink{omset}{\mathcal{W}}$ has \hyperlink{om2}{$(\omega_2)$}. Moreover, by \cite[Corollary 5.8]{compositionpaper} and \cite[Corollary 4.8]{testfunctioncharacterization}, one has that the following assertions are equivalent:
\begin{itemize}
\item[$(i)$] $\omega\in\hyperlink{omset}{\mathcal{W}}$ is quasianalytic,
\item[$(ii)$] $\Omega$ is quasianalytic in the sense of Definition \ref{matrixquasianalytic},
\item[$(iii)$] some/each $W^{(l)}$ satisfies (\ref{eq_Q}).
\end{itemize}

\medskip
Similarly, from \cite[Proposition 2]{borelmappingquasianalytic} (and in the same spirit as in \cite[Section 5]{compositionpaper}), one knows that $\Lambda^r_{[\omega]}=\Lambda^r_{[\Omega]}$ as locally convex spaces, too.\vspace{6pt}

Consequently, under the assumptions described above, we are able to apply the results from Section \ref{Borelmapfarawaymatrix} to the matrix $\mathcal{N}\equiv\Omega$, using the sequence $L$ from Proposition \ref{Thillieztheoremmatrixgeneralizationroum} lying above the matrix $\mathcal{M}\equiv\Sigma$ which is associated with a given (arbitrary) quasianalytic weight function $\sigma\in\hyperlink{omset}{\mathcal{W}}$.

\begin{theorem}\
\begin{itemize}
\item Let $\omega\in\hyperlink{omset}{\mathcal{W}}$ be a quasianalytic weight function satisfying $\hyperlink{om5}{(\omega_5)}$. Then, for any quasianalytic weight function $\sigma\in\hyperlink{omset}{\mathcal{W}}$, the set $j^{\infty}(\mathcal{E}^{0,r}_{\{\sigma\}}) \cap \Lambda^r_{(\omega)}$ is meager and Haar-null in $\Lambda^r_{(\omega)}$.
\item Let $\omega\in\hyperlink{omset}{\mathcal{W}}$ be a quasianalytic weight function satisfying \hyperlink{om2}{$(\omega_2)$} and $\liminf_{t\rightarrow+ \infty}\frac{\omega(t)}{t}=0$. Then, for any quasianalytic weight function $\sigma\in\hyperlink{omset}{\mathcal{W}}$, the set $\Lambda^r_{\{\omega\}} \setminus j^{\infty}(\mathcal{E}^{0,r}_{\{\sigma\}})$ is dense in $\Lambda^r_{\{\mathcal{\omega}\}}$ (and so $\Lambda^r_{\{\omega\}} \setminus j^{\infty}(\mathcal{E}^{0,r}_{(\sigma)})$ too).
\item Let $\omega\in\hyperlink{omset}{\mathcal{W}}$ be a quasianalytic weight function satisfying \hyperlink{om2}{$(\omega_2)$} and $\liminf_{t\rightarrow+ \infty}\frac{\omega(t)}{t}=0$ in the Roumieu resp. $\hyperlink{om5}{(\omega_5)}$ in the Beurling case. Then, for any quasianalytic weight function $\sigma\in\hyperlink{omset}{\mathcal{W}}$, the set $\Lambda^r_{[\omega]} \setminus j^{\infty}(\mathcal{E}^{0,r}_{\{\sigma\}})$ is lineable in $\Lambda^r_{[\omega]}$ (and so $\Lambda^r_{[\omega]} \setminus j^{\infty}(\mathcal{E}^{0,r}_{(\sigma)})$ too).
\end{itemize}

\end{theorem}

\bigskip

\textbf{Acknowledgement.} The authors wish to thank the referee for his useful comments and suggestions which have improved the presentation and the structure of this work.

C. Esser is supported by a F.R.S. -FNRSgrant; G. Schindl is supported by FWF-Project J3948-N35, as a part of which he is an external researcher at the Universidad de Valladolid (Spain) for the period October 2016 - September 2018.

\bibliographystyle{plain}
\bibliography{Bibliography}
\end{document}